\magnification=1200
\input amstex
\documentstyle{amsppt}
\hsize=5.25in
\vsize=7.25in
\nologo
\font\ninerm = cmr9

\TagsOnRight
\document
\def\sg{\sigma}
\def\LL{\Cal L}
\def\BB{\Cal B}
\def\AA{\Cal A}
\def\lm{\lambda}

\def\defeq{\hbox{$=$\hskip -10pt\raise 5.5pt\hbox{ def}}}
%

%

\centerline{\bf DOUBLE CENTRALIZING THEOREMS FOR}
\centerline{\bf THE ALTERNATING GROUPS}
\vskip 20pt
\centerline{\bf Amitai  Regev\footnote{Partially supported by ISF
Grant 6629 and by Minerva Grant No. 8441.}}
\centerline{\bf Department of Mathematics}\par
\centerline{\bf The Weizmann Institute of Science}\par
\centerline{\bf Rehovot 76100, Israel}\medskip
\centerline{\it E-mail:~~regev\@wisdom.weizmann.ac.il}\bigskip
%
%
\vskip 3cm
\centerline{\bf Abstract.}
\bigskip{\leftskip 30pt\rightskip 30pt\ninerm
Let $V^{\otimes n}$ be the $n$--fold tensor product of a vector space $V.$
Following I. Schur we consider the action of the 
symmetric group $S_n$ on $V^{\otimes n}$ 
by permuting coordinates. In the `super' ($\Bbb Z _2$ graded) case 
$V=V_0\oplus V_1,$ a $\pm$ sign is added [BR]. These actions 
give rise to the
corresponding Schur algebras S$(S_n,V).$ Here S$(S_n,V)$ is
compared with 
S$(A_n,V),$ the Schur algebra corresponding to  the 
alternating subgroup $A_n\subset S_n\,.$  While in the `classical' (signless) case these two 
Schur algebras are the same for $n$ large enough, it is proved that in the 
`super' case where $\dim V_0=\dim V_1\,,\;$
S$(A_n,V)$ is isomorphic to the crossed--product algebra
S$(A_n,V)\cong$ S$(S_n,V)\times\Bbb Z _2\,.$\par}
\newpage

\noindent {\bf\S 0. Introduction.}
\bigskip
Let $V$ be a finite dimensional vector space over the field $F=\Bbb C$ 
of the complex numbers 
(in fact, we shall only need the fact that 
$\sqrt {-1}\in F$), and let $V^{\otimes n}= V\otimes \cdots \otimes V$ 
$n$ times. The symmetric group $S_n$ acts on 
$V^{\otimes n}$ (say, from the left) by permuting coordinates.
This makes $V^{\otimes n}$ a left $FS_n$ module with the corresponding
 Schur algebra $\;End_{{FS_n}}(V^{\otimes n}).$ Here $FG$ denotes the 
group algebra of a group $G.$
Formally, that action is given by a multiplicative homomorphism 
$\varphi:\; S_n \longrightarrow End_{F}(V^{\otimes n})$ 
which extends linearly to an algebra homomorphism
$\varphi:\; FS_n \longrightarrow End_{F}(V^{\otimes n})\,.$ 
\medskip
Let $A_n\subset S_n$ denote the alternating group, with the 
corresponding Schur algebra $End_{{FA_n}}(V^{\otimes n}).$ The main purpose 
of this paper is to study and compare the pair of Schur algebras
$$End_{{FS_n}}(V^{\otimes n})\subseteq End_{{FA_n}}(V^{\otimes n}).$$
We mention first the following phenomena:
\bigskip 
\noindent {\bf Theorem 1 (see Remark 1.9).}
Let $\dim V=k$ and consider the above (signless!) action of 
$S_n$ on $V^{\otimes n}\,.$ 
If $k^2 \lvertneqq n$ then 
$$
\varphi(FA_n) = \varphi(FS_n),
\qquad\text{hence also}\qquad
End_{{FA_n}}(V^{\otimes n})=
End_{{FS_n}}(V^{\otimes n}). \tag 0.1
$$
\medskip
In a sense, Theorem 1 shows an anomaly: even though $|S_n|=2|A_n|\,,$ 
nevertheless $\varphi (FS_n)=\varphi (FA_n)\,,$ 
provided $k^2\lvertneqq n\,.$
As indicated below,
the incorporation of a $\pm$ sign to the above permutation action of $S_n$ 
seems to be natural and to remove that anomaly. Such a sign permutation 
action is related  to the representation theory of
Lie superalgebras [BR] [Se]. We now briefly explain that 
$S_n$ action, and this will allow us to formulate the main result of this
paper, which is Theorem 2 below.

\medskip
Let $V=V_0 \oplus V_1\,,\; \dim V_0=k\,,\; \dim V_1=l\,,\;$  
and let $S_n$ act on $V^{\otimes n}$ by permuting coordinates as before, but now,
together with a $\pm$ sign; that sign is obtained by considering 
the elements of $V_0$ as being
central, and the non-zero elements of $V_1$ as anti-commuting
among themself.
This is the sign-permutation action * of $S_n$ on $V^{\otimes n}$ 
[BR \S 1], [Se]. This action determines the algebra homomorphism
$\varphi ^*:\; FS_n \longrightarrow End_{F}(V^{\otimes n}).$ It  endows
$V^{\otimes n}$ with a new $FS_n$ module structure, which yields the
corresponding (new!) Schur algebras
$$
End_{FS_n}(V^{\otimes n})=End_{\varphi^*(FS_n)}(V^{\otimes n})
\defeq\BB _n\,,\qquad\text{and}  
$$
$$
End_{FA_n}(V^{\otimes n})=End_{\varphi^*(FA_n)}(V^{\otimes n})
\defeq\AA _n \,,
$$
and clearly $\BB _n\subseteq\AA _n.$
The main result of this paper is the following {\it crossed product theorem.}
\bigskip

\noindent {\bf Theorem 2 (The Cross Product Theorem. See Theorem 1.1).} 
Let $V=V_0 \oplus V_1$ with
$\dim V_0=\dim V_1$ and consider the above *-action of $S_n$ on 
$V^{\otimes n}.$
Then
$$
\dim \varphi ^* (FS_n) = 2 \dim \varphi ^* (FA_n)\quad \text{and}\quad  
\dim \AA _n=2 \dim \BB _n\,.
\tag 0.6
$$
\noindent   Moreover,  
there exists an algebra automorphism 
$\omega$  of order two,
$\omega \,: \BB_n \to \BB_n\,,$  such that $\Cal A _n$ is isomorphic to the 
crossed  product 
$\Cal A _n\cong \BB_n \times F[H]\cong \BB_n \times F[\Bbb Z_2],$ 
where $H=\{1,\omega\}.$
\bigskip
The case $\dim V_0\ne\dim V_1$ is considered in Section 3, and we prove
\medskip
\noindent {\bf Theorem 3 (Theorem 3.2).} 
Let $\dim V_0\ne\dim V_1$, then each of the above subalgebras of 
$End_F(V^{\otimes n})$ splits into a direct sum 
of two subalgebras as shown below, and the following relations hold between
these summands.
\medskip 
\noindent 1) 
$$
\;\varphi ^* ( FS_n)= P \oplus Q\quad\text{and} 
\quad\varphi ^* ( FA_n)= P' \oplus Q', \tag 0.7
$$ 
satisfying  $Q'=Q,\;\,$  $\,P'\subset P\,\;$ and $\;\dim P = 2 \dim P'.$
\medskip 
\noindent 2) 
$$
\;End_{{FS_n}}(V^{\otimes n})=\BB_n=\BB_P\oplus \BB_Q,\;
\;End_{{FA_n}}(V^{\otimes n})=\Cal A _n=\Cal A_P \oplus \Cal A _Q\,,\tag 0.8
$$
satisfying  $\BB_Q=\Cal A_Q$ and $\BB_P\subset\Cal A_P.$
\bigskip
\noindent At the moment, when $\dim V_0\ne\dim V_1$, the relations 
between $\BB _P$ and $\AA _P$ are not clear.
\bigskip
Section 2 contains some necessary preliminaries from the representation
theory of $S_n$ and $A_n$ -- which are applied in Section 3.
\bigskip
In Section 4 we first consider dimensions growth rates, 
as $n\to\infty,$  of the factors in (0.7) and (0.8). We show (Remark 4.1) 
that when  $\dim V_0\ne\dim V_1$,  
the above $P$-summands have much smaller dimensions than
the corresponding $Q$-summands. In particular, this implies that 
in that case,
$$
\lim_{n\longrightarrow \infty}
\frac{\dim\varphi ^* ( FA_n)}{\dim  \varphi ^*  ( FS_n)}= 
\lim_{n\longrightarrow \infty}
\frac{\dim\AA _n}{\dim \BB_n}= 
\lim_{n\longrightarrow \infty}
\frac{\dim End _{{FA_n}}(V^{\otimes n})}
{\dim End _{{FS_n}}(V^{\otimes n})}= 1\,.
\tag 0.9
$$

\noindent In the rest of Section 4 the case $\dim V_0=\dim V_1$ is revisited.
There, we study dimensions growth rates, 
as $n\to \infty\,,$ of certain factors that
appear in theorem 2. These factors are related to the self--conjugate
and to the non--self--conjugate partitions in the $(k,k)$ hook
$H(k,k;n)\,.$
It is shown that the asymptotics 
for the dimensions of these factors do exist, and can be
computed explicitly -- by Selberg integrals.
\bigskip

\noindent {\bf Remark.} An important step in the proof of Theorem 2 
(i.e. Theorem 1.1) is Proposition 
1.2. In its original proof, a considerable 
amount of the representation theory of $S_n$ and of $A_n$ 
was applied. Its present
elegant proof is due to G. Olshanski.
\bigskip \bigskip
\newpage

\noindent {\bf \S 1 The case $\dim V_0=\dim V_1\,$.}
\bigskip
\noindent Let $V=V_0 \oplus V_1$ be a $\Bbb Z _2$ 
graded vector space over a field 
$F$ of characteristic zero, and
consider the sign action * of $S_n$ on $V^{\otimes n}$ [BR] [Se].
Denote that action by $\varphi ^*\,.$ For example, let $n=2,\;\,
\sg=(1,2)\in S_2$ and let $y_1,y_2\in V$ be homogeneous, then 
$$
\varphi _{(1,2)}^*(y_1\otimes y_2)=
(-1)^{(\deg y_1)(\deg y_2)} (y_2\otimes y_1)\,.
$$
This yields the algebra homomorphism 
$$
\varphi ^*:\; FS_n \longrightarrow End_{_F}(V^{\otimes n}) 
$$
and defines the centralizer $\BB _n=End_{_{FS_n}}(V^{\otimes n})\;\;$ 
of $\;\;\varphi ^*(FS_n)\;\;$ in $\;\;End_{_F}(V^{\otimes n}).$
Obviously, for the alternating group $A_n \subset S_n$ we have
$$
\varphi ^*(FA_n)\subseteq \varphi ^*(FS_n)\qquad\text{and}\qquad
End_{_{FS_n}}(V^{\otimes n})\subseteq End_{_{FA_n}}(V^{\otimes n}),
$$
and we denote  $\AA _n=End_{_{FA_n}}(V^{\otimes n})\,.$
Thus, $\;\AA _n\;$ is the centralizer of 
$\,\varphi ^*(FA_n)\;$ in 

\noindent $End_{_F}(V^{\otimes n}).$
We are interested in the comparison between $\BB _n$ and $\Cal A _n\,$
and also between the images $\varphi ^* (FA_n)$ and 
$\varphi ^* (FS_n)\,.$ 
In this section we study the case when $\dim V_0=\dim V_1.$ 
The main result of this paper is the following theorem. 
\bigskip
\noindent {\bf Theorem 1.1.}$\;$ Let $V=V_0 \oplus V_1,\;\,
\dim V_0=\dim V_1 \,$ and consider the above sign permutation action *
of $S_n$ on $V^{\otimes n}\,.$
Then there is a subgroup $H$ of the automorphisms of 
$\BB _n=End_{FS_n}(V^{\otimes n}) \,,\;\,
H\cong \Bbb Z _2\,,$ such that the algebra 
$\Cal A _n=End_{FA_n}(V^{\otimes n})$ 
is isomorphic to  a crossed product of the algebra
$\BB_n=End_{FS_n}(V^{\otimes n})$ with the group $H\,.$
In particular, 
$$
\dim \left ( End_{_{FA_n}}(V^{\otimes n})\right )=
2 \dim \left ( End_{_{FS_n}}(V^{\otimes n})\right )\quad\text{i.e.}\quad
\dim  \Cal A _n= 2\dim \BB_n\,.            \tag 1.1.1
$$
\noindent {\bf Remark.} By applying the representation theory of $A_n$ 
we show in the next section that in the situation of Theorem 1.1, also
$$
 \dim \varphi ^* (FS_n) = 2 \dim \varphi ^* (FA_n).    \tag 1.1.2
$$
\medskip
The Key for proving Theorem 1.1 is 
\bigskip
\noindent {\bf Proposition 1.2.} Let 
$V=V_0\oplus V_1,\;\;\dim V_0=\dim V_1\,,$ 
len $n \ge 2$ and let $S_n $ 
act with the sign permutation action *
on $V^{\otimes n}$ and let $\AA _n= End_{FA_n}(V^{\otimes n})$ and
$\BB _n = End_{FS_n}(V^{\otimes n})\,,$ then 
$$
\dim \AA _n = 2\dim \BB _n\,. \tag 1.2.1
$$
The following proof of Proposition 1.2 was suggested to us by G. Olshanski. 
\bigskip
Consider first the general 
$(k,l)$ case: let 
$\dim V_0=k,\;\dim V_1=l,\; V=V_0 \oplus V_1,$ and
consider the action * of $S_n$ on $V^{\otimes n}$ [BR] [Se].
\medskip

\noindent {\bf Definition 1.3.} 
Denote by $I _n^+(V)$ the subspace of $V^{\otimes n}$
of the symmetric tensors, and by $I _n^-(V)$ the subspace of $V^{\otimes n}$
of the anti-symmetric tensors. Thus 
$$
I _n^+(V)=\{w\in V^{\otimes n}\;|\; \varphi _{\sigma} ^*(w)=w
\quad\text{for all}\quad \sigma\in S_n\}\tag 1.3.1
$$
and
$$
I _n^-(V)=\{w\in V^{\otimes n}\;|\; \varphi _{\sigma} ^*(w)=\text{sgn}
(\sigma)w
\quad\text{for all}\quad \sigma\in S_n\}\,.\tag 1.3.2
$$
Finally, denote
 by $\LL _n(V)$ the subspace of $V^{\otimes n}$ of the 
$A_n$ -- invariant vectors:
$$
\LL _n(V)=\{w\in V^{\otimes n}\;|\; \varphi _{\sigma} ^*(w)=w
\quad\text{for all}\quad \sigma\in A_n\}\,.\tag 1.3.3
$$
\medskip
\noindent {\bf Lemma 1.4.} We have 
$$
\LL _n(V)=I _n^+(V)\oplus I _n^-(V)\,.\tag 1.4.1
$$
\medskip
\noindent {\bf Proof.} The inclusion $\supseteq$ is obvious. The opposite 
inclusion is obtained from
$$
w=\frac{1}{2} (w+\varphi _{(1,2)} ^*w)+ \frac{1}{2}  (w-\varphi _{(1,2)} ^*w),
\tag 1.4.2
$$
where $w\in \LL _n(V)\,.$
\bigskip
\noindent {\bf Remark 1.5.} It is well known that if 
$l=\dim V_1=0$ and $k=\dim V_0\lneqq n$ then $I _n^-(V)=0\,.$
\bigskip
We show next that when $\dim V_0=\dim V_1,\;\,I _n^+(V)\cong I _n^-(V)$
(Lemma 1.8). In fact, we produce an explicit such isomorphism
$T\,.$ To do so, 
denote 
$$
E=End_F(V)
$$ 
and note it is $\Bbb Z_2$ graded.
\bigskip
\noindent {\bf The identification 
$E^{\otimes n}\equiv End_F(V^{\otimes n})\,.$ }
We use standard sign convention for the 
canonical isomorphism $X\otimes Y \to Y\otimes X$, that is, for homogeneous 
$x,y,$
$$
x\otimes y \to (-1)^{\deg\,x\deg\,y}y\otimes x\,.
$$
This induces  a corresponding identification 
$E^{\otimes n}\equiv End_F(V^{\otimes n})\,.$ 
\smallskip
\noindent  For example, let $n=2\,,\;\,f_1\,,\,f_2\in E$ homogeneous
and $v_1\,,\,v_2\in V$ homogeneous, then
$$
(f_1\otimes f_2)(v_1\otimes v_2)=(-1)^{(\deg(f_2))(\deg(v_1))}
(f_1(v_1)\otimes f_2(v_2))\,.
$$
When $V$ is replaced by $E\,,$ the action $\varphi ^*$ of $S_n$ on
$V^{\otimes n}$ is replaced by the corresponding 
action $\eta ^*$ of $S_n$ on
$E^{\otimes n}\,.$
These two actions are related by formula (1.6.1) below.
\medskip
\noindent {\bf Lemma 1.6. [Se, Lemma 1].} 
Let $L\in End_F(V^{\otimes n})\,\;\,$ 
and $\sg\in S_n\,,$ then 
$$
\eta _{\sg}^*L=\varphi _{\sg}^*\circ L\circ \varphi _{\sg ^{-1}}^*\,.
\tag 1.6.1
$$
This obviously implies
\medskip
\noindent {\bf Corollary 1.7.} Let $E=End_F(V),$
then 
$$
I _n^+(E)=End_{FS_n}(V^{\otimes n})=\BB _n \tag 1.7.1
$$
and
$$
\LL _n(E)=End_{FA_n}(V^{\otimes n})=\AA _n\,.\tag 1.7.2
$$
\bigskip
\noindent {\bf Lemma 1.8.} Let  
$\dim V_0=\dim V_1\,.$ Then there exist $T\in \AA _n$ which 
induces an isomorphism between 
$I _n^+(V)$ and  $I _n^-(V)\,:$ $T( I _n^+(V))=I _n^-(V)\;$
and $\;T(I _n^-(V))=I _n^+(V)\,.$ Moreover, 
$$
T^2= (-1)^{\binom n2}\cdot
I\quad(\text{$I$ is the identity map}),\tag 1.8.1
$$
and 
$$
T\circ\varphi_{(1,2)}^*=- \varphi_{(1,2)}^*\circ T\,.\tag 1.8.2
$$
\medskip
\noindent {\bf Proof.} 
First, construct $\tau\,:$ Let
$t_1,\cdots , t_k\in V_0\;, u_1,\cdots , u_k\in V_1$ be bases, and let
$\tau\in  E$ be given by 
$\tau t_i=u_i$ and $\tau u_i=t_i\;, i=1,\cdots , k\,.$
Note that $\tau \in E_1\,,$ i.e. $\deg\tau=1\,,$ 
where $E=End_F(V)=E_0\oplus E_1\,.$

Now define $T=\tau ^{\otimes n}.$ For example, if $n=2$ and 
$y_1,y_2 \in V$ are homogeneous, 
$$
T(y_1\otimes y_2)=  (\tau\otimes \tau) (y_1\otimes y_2)=
(-1)^{\deg y_1}\cdot \tau y_1\otimes \tau y_2\,.
$$
Since $\tau \in E_1\,,$ it follows from the definition of $\eta _{\sg}^*$
that for any $\sg\in S_n,\;\,\eta _{\sg}^*T=\text{sgn}(\sg)T\,.$ 
By (1.6.1) it follows that 
$$
\varphi _{\sg}^*\circ T=\text{sgn}(\sg)T\circ\varphi _{\sg}^*\,,\qquad
\sg\in S_n\,,
$$
which clearly implies that $T\in \AA _n$ and is an isomorphism between 
$I _n^+(V)$ and $I _n^-(V)\,.$  
Property (1.8.2) is now obvious, while (1.8.1) follows from the definition
of $T$ and from the fact that if $y\in V$ is homogeneous, then 
$$
\deg (y) + \deg (\tau (y)) =1\,.
$$ \hfill \qed
\bigskip
We can now complete
\medskip
\noindent {\bf The proof of Proposition 1.2.} 
Apply (1.7.2), then apply (1.4.1)  with $E$ replacing $V\,,$ to deduce that 
$\AA _n =I _n^+(E)\oplus I _n^-(E)\,.$ By (1.7.1) $\BB _n =I _n^+(E),$
and the proof follows from Lemma 1.8, again with $E$ replacing $V\,.$ 
\hfill \qed
%
\bigskip
\noindent {\bf Remark 1.9.} Let $\dim V_1=0\,,$ so that the sign permutation
action $\varphi ^*$ is the classical  permutation action
$\varphi$ of $S_n$ on $V^{\otimes n}$, as considered by Schur. Let 
$\dim V=k$, hence $\dim E=k^2,$ and let $k^2<n\,.$ By Remark 1.5 
and by (1.4.1), with $E$
replacing $V\,,$ $\LL _n(E)=I _n(E),$ and by (1.7.1) and (1.7.2), it follows 
that 
$$
\AA _n =\BB _n\,.\tag 1.9.1
$$
 Also, by double centralizing, this implies that 
in that case, 
$$
\varphi (FS_n)=\varphi (FA_n)\,.\tag 1.9.2
$$
As was remarked in the Introduction,
both (1.9.1) and (1.9.2) may be considered as a sort of anomaly, since one 
would expect $\dim\varphi (FS_n)$ to be twice as large as 
$\dim \varphi (FA_n)\,,$ and similarly for $\dim \AA _n$ versus 
$\dim \BB _n\,.$ Note that 
this anomaly disappears when 
the ordinary case $\varphi$  is replaced by the super case $\varphi ^*\,,$ 
where $\dim V_0=\dim V_1\,.$ 
\bigskip
We  now complete the proof of Theorem 1.1.
\bigskip
\noindent {\bf Lemma  1.10.} Let  
$\dim V_0=\dim V_1\,$ and $T\in \AA _n$ as in Lemma 1.8, with $T$ satisfying
(1.8.1) and (1.8.2). Then 
$$
\AA _n = \BB _n \oplus T \BB _n \,,\tag 1.10.1
$$
$T\BB _n T=\BB _n\,,\;$ and the map 
$$
\omega : \BB_n \longrightarrow \BB_n\,,\quad\quad \quad
\omega (b)= \varepsilon\cdot  TbT=
(-1)^{\binom n2}\cdot  TbT,
\quad (\varepsilon=(-1)^{\binom n2}) \tag 1.10.2 
$$
is an automorphism of order two of $\BB _n\,.$ 
\smallskip
\noindent We denote $H=\{1,\omega\},$ which is clearly a group isomorphic
to $\Bbb Z_2\,.$
\bigskip
\noindent {\bf Proof.} It easily follows from (1.8.1) and (1.8.2) that 
$\BB _n \cap T \BB _n=0\,,$ therefore 
$\BB _n \oplus T \BB _n\cong \BB _n + T \BB _n\subseteq \AA _n\,,$
and equality follows by the equality of the dimensions, as implied 
by Proposition 1.2  and the obvious equality 
$\dim \BB _n=\dim T\BB _n\,.$ This proves (1.10.1). Also,  (1.10.2)
(i.e. the fact that $\omega$ is of order two)
easily follows from (1.8.1). 
\bigskip
\noindent {\bf 1.11 . Crossed products.} 
Let $G$ be a finite group of automorphisms of an $F$-algebra 
$R:\;\, g\in G\,;r\in R,\; g:r \to g(r).$ We denote the crossed product
of $R$ and $G$ by $R\times FG,$ where $FG$ is the group algebra.
Recall that as an $F$  vector space, 
$R\times FG=R\otimes _F FG.$ Multiplication is defined by linearity and
by $(r_1\otimes  g_1)(r_2\otimes  g_2)=r_1\cdot g_1(r_2)\otimes g_1g_2\,.$

Let now $\omega$ be given by (1.10.2), let $H=\{1,\omega\}$
and let $\BB _n\times FH$ be the corresponding crossed product. Thus,
$B_n\times FH =B_n\otimes _FFH,$ and multiplication is given 
as follows. Let $b_1,\,b_2\in B_n,\; s\in FH,$ then
$$
(b_1\otimes  1)(b_2\otimes  s)=b_1b_2\otimes s\;\, \text{and} \;\,
(b_1\otimes  \omega )(b_2\otimes   s)=b_1\omega(b_2)\otimes  \omega s=
b_1Tb_2T\otimes  \omega s\,.
$$
\medskip
We now reformulate -- and prove -- Theorem 1.1. This is
\medskip
\noindent {\bf Theorem 1.12 (1.1').} Let $\dim V_0=\dim V_1$ and let  
$\AA _n=End_{FA_n}(V^{\otimes n})$
and $\BB _n=End_{FS_n}(V^{\otimes n})$ (with respect to 
$\varphi ^*\,).$ Let  $\omega$ be
given by (1.10.2), $H=\{1,\omega\},$  and $\BB _n\times FH$ the corresponding
crossed product.  Then there is an algebra isomorphism
$$
\AA _n \cong \BB _n\times FH\;\;(\,\cong 
\BB _n\times F\Bbb Z _2 \,)\, . \tag 1.12.1
$$
Explicitly, the isomorphism is given by
$$
\zeta: \; \BB _n\times FH \longrightarrow \AA _n \,,
$$
where $\zeta$ is given by linearity and by 
$$
\zeta(b\otimes 1)=b\qquad \text{and}\qquad \zeta(b\otimes \omega)=
\sqrt{\varepsilon}\cdot bT,
\tag 1.12.2
$$
where, again, $\varepsilon=(-1)^{\binom n2}\,.$
\medskip \medskip

\noindent {\bf Proof.} It easily follows from Lemma 1.10 that 
the above 
$\;\zeta$ is a vector space isomorphism.
It is easy to see that it is also an algebra homomorphism. For example,
$$
\multline
\zeta ((b_1 \otimes \omega)(b_2 \otimes 1))= 
\zeta  ((b_1 \omega (b_2)\otimes \omega))=\\
\zeta  ((b_1\varepsilon  Tb_2T\otimes \omega))=
\varepsilon \sqrt{\varepsilon}
b_1Tb_2TT=\sqrt{\varepsilon} b_1Tb_2,  
\endmultline
$$
while, also,  
$$
\zeta (b_1 \otimes \omega)\zeta (b_2 \otimes 1)=\sqrt{\varepsilon} b_1Tb_2.
$$

Similarly,
$$
\multline
\zeta ((b_1 \otimes \omega)(b_2 \otimes\omega ))= 
\zeta  ((b_1 \omega (b_2)\otimes 1))=\\
\zeta  ((b_1\varepsilon  Tb_2T\otimes 1))=
\varepsilon b_1Tb_2T, 
\endmultline
$$
while, also,
$$
\zeta (b_1 \otimes\omega )\zeta(b_2 \otimes \omega)=
(\sqrt{\varepsilon})^2 b_1 Tb_2T.
$$
Similarly for the other (trivial) two cases 
$(b_1 \otimes 1)((b_2 \otimes \omega)$ and 
$(b_1 \otimes 1)((b_2 \otimes 1).$
\hfill \qed
\bigskip
\bigskip

\noindent {\bf \S 2. Some $S_n$ and $A_n$ Representations.}
\bigskip
\noindent We summarize some facts (and notations) about 
$S_n$ and $A_n$ representations -- that will be needed later 
[B] [J] [JK] [M] [Sa].
\smallskip
\noindent For 
$\lambda =(\lambda _1,\cdots , \lambda _r )  \vdash n$ (i.e.
$\lambda$ is a partition of $n$) denote $l(\lambda )=r$ if
$\lambda _r \ne 0\,.$ Also,
$\lambda '$ is the conjugate partition of $\lambda\,.$
\medskip
\noindent {\bf 2.1.} The group algebra $\;FS_n$ admits the following 
decomposition:
$$
FS_n=\bigoplus _ {\lambda \vdash n}I_ {\lambda },\tag 2.1.1
$$
where each $I_ {\lambda }$ is a minimal two-sided ideal in $FS_n\,,$
and is isomorphic to certain matrix algebra
$M_r(F)$ ($r\times r$ matrices over $F$). Moreover,
there is an explicit correspondence $\lambda \longleftrightarrow I_ {\lambda },$
due to Frobenius and to A. Young,  for which 
$I_ {\lambda} \cong M_{f^ {\lambda}}(F),$ where 
$f^ {\lambda}$ is the number of standard Young tableaux of
shape ${\lambda}.$ Also, 
$I_ {\lambda }
\cong J_{\lambda}^{\oplus {f^ {\lambda}}},$ where 
$J_{\lambda}$ is a  left ideal, irreducible as an $FS_n$ left module.

\noindent Moreover, if $\lambda\ne \mu \vdash n$ then 
$J_{\lambda} \not\cong J_{\mu}$ as $FS_n$ left modules.
\medskip

\noindent {\bf 2.2.} Consider $FS_n$ as a left $FA_n$ module.

\noindent {\bf 2.2.1.} If $\lambda\ne \lambda '$ then, as  $FA_n$ modules,
$J_{\lambda} \cong J_{\lambda '} $ and both remain irreducible.

\noindent{\bf 2.2.2.} If $\lambda= \lambda '$ then, as  $FA_n$ modules,
$J_{\lambda} =J_{\lambda}^+ \oplus J_{\lambda}^- ,\;J_{\lambda}^+ $
and $J_{\lambda}^- $ are $FA_n$ irreducible, and
$J_{\lambda}^- =(1,2) J_{\lambda}^+$ for the transposition 
$(1,2)\in S_n\,.$ In particular,
$\dim  J_{\lambda}^+ =\dim  J_{\lambda}^-=\frac{1}{2}\dim  J_{\lambda}.$
\bigskip
\noindent {\bf 2.3.} Note  that 2.1
is the decomposition of the $S_n$ regular representation into 
{\it isotypic} components: each $I_{\lambda}$ is the sum of all the 
irreducible $FS_n$ submodules of $FS_n$  which are isomorphic to a 
given irreducible 
$FS_n$ module. Thus, if $M\subseteq FS_n $ is an $FS_n $
left submodule and $M \cong J_{\lambda}$ 
then $M\subseteq I_{\lambda}$
($J_{\lambda},\;I_{\lambda} $ as in 2.1).
\smallskip
\noindent {\bf 2.4.} Let `$\mu <\lambda $' denote the left lexicographic order
on the partitions of $n$. It follows from 2.2 that the isotypic decomposition 
of $FS_n$ as a left $FA_n$ module is given by:
$$
FS_n=\left [\bigoplus_{\lambda \vdash n; \lambda > \lambda '}
(I_{\lambda} \oplus I_{\lambda '})\right ]\oplus 
\left [\bigoplus_{\lambda\vdash n ; \lambda =\lambda '}
(I^+_{\lambda}\oplus I^-_{\lambda})\right ], \tag 2.4.1
$$
where for $\lambda \ne\lambda ',\;\,I_{\lambda}\oplus I_{\lambda '}\cong 
J_{\lambda}^{\oplus 2f^{\lambda}},\;\,$ and for   
$\lambda = \lambda ' ,\;\,I_{\lambda}=I_{\lambda} ^+ 
\oplus I_{\lambda} ^-,\;\,$
$\;\,I_{\lambda} ^+\cong\left (J_{\lambda} ^+\right )^
{\oplus f^{\lambda}}\;\,$ and $\;\,
I_{\lambda} ^-\cong\left (J_{\lambda} ^-\right )^
{\oplus f^{\lambda}}\,.$ 
\smallskip
\noindent Similarly, it follows that the isotypic decomposition 
of the regular representation of $A_n$ is 
$$
FA_n=\left [\bigoplus_{\lambda \vdash n; \lambda > \lambda '}
\bar I_{\lambda}\right ]\oplus \left [\bigoplus_{\lambda \vdash n; 
\lambda =\lambda '}
(\bar I_{\lambda}^{+}\oplus \bar I_{\lambda}^{-} )\right ].\tag 2.4.2 
$$
Here, for $\lambda\ne\lambda '\,, \;\bar I_{\lambda}$ 
is the sum of all the irreducible left  
$FA_n$ submodules 
of $FA_n$ which are isomorphic to $J_{\lambda}\,.$
Similarly for $\bar I_{\lambda}^{+}\;$ and for $\;\bar  I_{\lambda}^{-} .$
Recall that for any group $G,$ the decomposition of the  regular
representation into isotypic components -- coincide with the decomposition
of $FG$ into simple two-sided ideals. 
Comparing the isotypic decompositions of $FS_n$ and of $FA_n$ as $FA_n$
modules, we deduce 
\medskip
\noindent {\bf Lemma 2.5} 
The above (2.4.1) gives the decomposition of $FS_n$ into simple 
two-sided ideals,
while  (2.4.2) gives that of $FA_n$ into simple two-sided ideals.
These decompositions satisfy
\smallskip
\noindent a) If $\lambda \ne \lambda '$ then 
$ \bar I_{\lambda}\subseteq I_{\lambda}\oplus I_{\lambda '}$ 
and as $FA_n$ modules, 
$\bar I_{\lambda}\cong I_{\lambda}\cong I_{\lambda '}.$
In particular, $\dim \bar I_{\lambda}=\dim I_{\lambda}=\dim I_{\lambda '}.$ 
\smallskip
\noindent b) If $\lambda =\lambda '$ then 
$\bar I_{\lambda}^{+}\oplus \bar I_{\lambda}^{-}\subseteq I_{\lambda}$
and  $ \dim\bar  I_{\lambda}^{+}= \dim\bar  I_{\lambda}^{-}=
\frac{1}{4}\dim I_{\lambda}$.
\bigskip

\noindent  {\bf 2.6. The `hook' theorem for $\varphi ^*.$} The following 
theorem is a generalization of a classical theorem of Schur and Weyl [W].
\medskip
\noindent {\bf Theorem [BR, 3.20].} Let
$$
H(k,l;n)=\{(\lambda _1,\lambda _2, \cdots)\vdash n\;\;|
\;\;\lambda _{k+1}\le l\}\tag 2.6.1
$$
denote the partitions of $n$ whose diagrams are contained 
in the $k,l$ hook. As in the previous section, let
$V=V_0\oplus V_1\,,\;\,\dim V_0=k$ and $\dim V_1=l\,,$ and 
let  $S_n$ act on $V^{\otimes n}$ 
by the sign permutation action, denoted by $\varphi ^*.$ Then 
$$
\varphi ^*(I_{\lambda})\;\cong\;\left\{
\aligned &I_{\lambda}\quad\text{if}\quad\lambda\in H(k,l;n)\\
&0\quad\text{otherwise}
\endaligned
\right .
\tag 2.6.2
$$
and, moreover, 
$$
\;\varphi ^*(FS_n)=\bigoplus _{\lambda \in H(k,l;n)}\varphi ^*(I_{\lambda})
\cong \bigoplus _{\lambda \in H(k,l;n)}I_{\lambda}\,.
\tag 2.6.3
$$ 

\noindent  {\bf 2.7. The decomposition of $W=  V^{\otimes n}$
and of $\BB _n\,.$}    
The above $S_n$ action $\varphi ^*$  on $W=V^{\otimes n}$
makes $W$ into a left $FS_n$ module, and the decomposition (2.6.3) 
implies a corresponding decomposition of $W\,:$ 
$$
V^{\otimes n} =\bigoplus _{\lambda \in H(k,l;n)}  W_{\lambda},
\tag 2.7.1
$$ 
where $W_{\lambda }$ is the $FS_n$ isotypic component of $V^{\otimes n}$
corresponding to $I_{\lambda}.$ Thus 
$$
I_{\lambda} W_{\lambda }\;\defeq \;
\varphi ^*(I_{\lambda}) W_{\lambda }\;=\;W_{\lambda }
\tag 2.7.2
$$
for $\lambda\in H(k,l;n)\,.$
Also,
$$
I_{\lambda} W_{\mu}=0\quad\text{if}\quad \lambda \ne \mu\vdash n\,.
\tag 2.7.3
$$
Similarly,
$$
\BB_n=End_{_{FS_n}}(V^{\otimes n})=
\bigoplus_{\lambda \in H(k,l;n)} \BB_{\lambda}\,, \tag 2.7.4 
$$
where 
$$
\BB_{\lambda}=End_{I_{\lambda}}(W_{\lambda })\defeq
End_{\varphi ^*(I_{\lambda})}(W_{\lambda })\tag 2.7.5
$$
(and it is well known that 
$End_F(W_{\lambda })\cong \varphi ^*(I_{\lambda})\otimes _F\BB _{\lambda}$).

\bigskip
\bigskip
\noindent {\bf \S 3. The case $\dim V_0 \ne \dim V_1\,$.}
\bigskip
\noindent Let $\,V=V_0 \oplus V_1$. In this section we compare 
$\varphi ^*( FS_n)\;\,$ with  $\;\,\varphi ^*( FA_n)\,.$ 
If $\,\dim V_0=\dim V_1\,,$ we show that 
$\dim\varphi ^*( FS_n)=2\dim\varphi ^*( FA_n)\,.$ Recall also that in that 
case, $\dim \AA _n=2 \dim \BB _n$ (Proposition 1.2).
In the general case, when  $\,\dim V_0\ne\dim V_1\,,$ the picture is more 
involved, and depends on the following decomposition of $H(k,l;n)\,.$
\bigskip

\noindent {\bf Definition 3.1.} Denote
$H_0(k,l;n)=\{\lambda \;|\;\lambda,  \lambda '\in H(k,l;n)\}$ and let 
$H_1(k,l;n)=H(k,l;n)\backslash H_0(k,l;n)\,,$ so that
$$
H(k,l;n)= H_0(k,l;n)\cup H_1(k,l;n)\,,\quad\text{a disjoint union.}  \tag 3.1.1
$$
\medskip
\noindent If $\lambda \in H_1(k,l;n)$ then, by definition, 
$\lambda '\not\in H(k,l;n)\,.$
Assume, without loss of generality, that  $l\leq k\,.$ If $l=k$ then 
$H(k,k;n)= H_0(k,k;n)$ and $H_1(k,k;n)=\varnothing.$ If $l<k,$ it is easy to see
that $H_0(k,l;n)=\{\lambda\in H(k,l;n)\; |\; \lambda _{l+1}\leq k\}$
(and $H_0(k,l;n)$ is close to $H(l,l;n)$ in an obvious sense).
\bigskip
Recall that $\BB_n=End_{FS_n}(V^{\otimes n})$ and 
$\AA_n=End_{FA_n}(V^{\otimes n})\,.$
\medskip
\noindent The disjoint union decomposition (3.1.1) 
implies the following ``$P-Q$'' direct sum 
decompositions. 
\bigskip
\noindent {\bf Theorem 3.2.} Let 
$\,\dim V_0\ne\dim V_1\,.$
The decomposition (3.1.1) implies that 
\medskip
\noindent (a) $\;\varphi ^* ( FS_n)= P \oplus Q ,\quad
\varphi ^* ( FA_n)= P' \oplus Q',\;$ 
satisfying $\,Q'=Q\,,\quad P'\subset P\,,$ and $\dim P=2\dim P'\,.$
\medskip 
\noindent (b) Also corresponding to (3.1.1) we have 
$\AA_n=\AA_P\oplus \AA_Q$ and $\BB_n=\BB_P\oplus \BB_Q\,,$
satisfying $\BB_P\subset \AA_P$ and $\BB_Q=\AA_Q\,.$
\smallskip

\noindent We remark that beside the fact that 
$\dim \BB_P < \dim \AA_P\,,$ at the moment, the relations 
between $\AA_P$ and $\BB_P$ are not clear. 

\noindent We remark also that in 
Section 4 we show that if $\,\dim V_0\ne\dim V_1$ then, asymptotically,
both $\dim\varphi ^* ( FS_n)\cong \dim Q'=\dim Q\cong\dim 
\varphi ^* ( FA_n)$ 
and $\dim \AA_n\cong\dim \AA_Q=\dim \BB_Q\cong\dim \BB_n\,.$

\bigskip

We proceed with the details. 
\medskip 
\noindent By (2.6.3) and (3.1.1),
$$
\varphi ^* ( FS_n)= P \oplus Q\, ,                            \tag 3.2.1
$$
$$
\multline
\text{where}\quad P=\bigoplus_{\lambda \in H_0(k,l;n)}\varphi ^*(I_{\lambda})
\cong\bigoplus_{\lambda \in H_0(k,l;n) } I_{\lambda} ,\\
\text{and}\quad Q=\bigoplus_{\lambda \in H_1(k,l;n)}\varphi ^*(I_{\lambda})
\cong\bigoplus_{\lambda \in H_1(k,l;n)} I_{\lambda}\,.
\endmultline
$$
By pairing together $\lambda ',\lambda \in H_0(k,l;n)$ (and choosing
the notation such that $\lambda '\leq\lambda $ in the left lexicographic
order), the term $P$ can be rewritten as
$$
P\cong \bigoplus_{\lambda \in H_0(k,l;n)} I_{\lambda}=
\left [\bigoplus_{\lambda '< \lambda \in H_0(k,l;n) }
(I_{\lambda} \oplus I_{\lambda '})\right ]\oplus 
\left [\bigoplus_{\lambda '=\lambda \in H_0(k,l;n)}
I_{\lambda} \right ].\tag 3.2.2
$$
Similarly,
$$\varphi ^*(FA_n )=P'\oplus Q' \,,\tag 3.2.3
$$ 
with $P'$ corresponding to $H_0(k,l;n)$ and $Q'$ -- to $H_1(k,l;n)\,.$
By Lemma 2.5, 
$$
\multline 
P'\cong 
\left [\bigoplus_{\lambda '< \lambda \in H_0(k,l;n) }
I_{\lambda} \right ]\oplus 
\left [\bigoplus_{\lambda '=\lambda \in H_0(k,l;n)}
(\bar I_{\lambda}^+\oplus \bar I_{\lambda}^-) \right ]\\
\text{and}\qquad Q'\cong\bigoplus_{\lambda \in H_1(k,l;n)}\bar I_{\lambda}\,.
\endmultline
\tag 3.2.4
$$
Similarly, (2.7.1) and (3.1.1) imply that
$$
V^{\otimes n}=W_P \oplus W_Q \,,                            \tag 3.2.5
$$
$$
\text{where}\qquad  W_P=\bigoplus_{\lambda \in H_0(k,l;n) } W_{\lambda} ,\quad 
\quad\text{and}\quad\quad W_Q=\bigoplus_{\lambda \in H_1(k,l;n)} W_{\lambda}\,.
$$
Also here we can rearrange terms and rewrite
$$
W_P= \bigoplus_{\lambda \in H_0(k,l;n)} W_{\lambda}=
\left [\bigoplus_{\lambda '< \lambda \in H_0(k,l;n) }
(W_{\lambda} \oplus W_{\lambda '})\right ]\oplus 
\left [\bigoplus_{\lambda '=\lambda \in H_0(k,l;n)}
W_{\lambda} \right ].\tag 3.2.6
$$

%
\noindent {\bf Lemma 3.3.} Let $\lambda\in H_1(k,l;n),$ and refer to 
(2.4.1) and (2.4.2). Then $\varphi ^*(I_{\lambda '})=0$ and
$\varphi ^*(\bar I_{\lambda })=\varphi ^*(I_{\lambda })\,.$
Thus, by (3.2.1) and (3.2.3), $Q=Q'\,.$
\bigskip
\noindent {\bf Proof.} Since $\lambda '\not\in H(k,l;n),\;\,
\varphi ^*(I_{\lambda '})=0\,.$ Also, $\lambda \ne  \lambda '$
therefore $\bar I_{\lambda }\cong I_{\lambda }$ and
$\bar I_{\lambda }\subset I_{\lambda }\oplus I_{\lambda '}\,,$ and hence 
$\varphi ^*(\bar I_{\lambda })\subseteq \varphi ^*( I_{\lambda })\,.$
To prove equality, it clearly suffices to show that 
$\varphi ^*(\bar I_{\lambda })\ne 0\,.$
\smallskip
\noindent Assume $\varphi ^*(\bar I_{\lambda })=0\,,$ so
$\varphi ^*(\bar I_{\lambda }) W_{\lambda }=0\,.$   
If $\mu\ne\lambda,$ then 
\smallskip
\noindent (1) if $\mu\ne\mu '\,,\quad\varphi ^*(\bar I_{\mu }) 
W_{\lambda }=0$  since 
$\varphi ^*(\bar I_{\mu })\subseteq \varphi ^*( I_{\mu }\oplus  I_{\mu ' });$
(the case $\mu '=\lm ' $ also follows, since $\varphi ^* (I_{\lambda '})=0$)
\smallskip
\noindent (2) if $\mu=\mu '\,,\quad 
\varphi ^*(\bar I_{\mu }^+ \oplus    \bar I_{\mu }^- )W_{\lambda }=0$
since $\varphi ^*(\bar I_{\mu }^+\oplus\bar I_{\mu }^-)
\subset \varphi ^*( I_{\mu })\,.$
\smallskip
\noindent We deduce that the assumption that 
$\varphi ^*(\bar I_{\lambda })=0$ implies that 
$\varphi ^*(FA_n)W_{\lambda }=0\,,$ which is a contradiction, since
$1\in FA_n$ and $\varphi ^*(1)W_{\lambda }=W_{\lambda }$ for any
$\lambda\in H(k,l;n)\,.$ \hfill\qed
\bigskip
\noindent {\bf Proposition 3.4.} Let $P,\;P'$ as in (3.2.2) and (3.2.3). Then
$P'\subset P$ and $\dim P=2\dim P'\,.$
\bigskip
\noindent {\bf Proof.} The inclusion $P'\subset P$ follows easily from 
Lemma 2.5. 

\noindent Let $\lm\in H_0(k,l;n)\,.$ If $\lm\ne\lm '$ then $P$ contains the
summand $\cong I_{\lm}\oplus I_{\lm '}\;\,(\dim  I_{\lm '}=\dim  I_{\lm })$
while $P'$ contains the summand $\cong \bar I_{\lm}\cong I_{\lm}\,.$

\noindent If $\lm=\lm '$ then $P$ contains the summand $\cong I_{\lm}$
while $P'$ contains the summand $\cong \bar I_{\lm}^+\oplus I_{\lm}^-$
and by Lemma 2.5.b, $\dim \bar I_{\lm}^+\oplus  I_{\lm}^-=
\frac{1}{2} \dim I_{\lm}\,.$ This obviously implies the proof. 
\hfill \qed
\bigskip
\noindent {\bf Note} that Lemma 3.3 and Proposition 3.4 prove part (a) of
Theorem 3.2
\bigskip
\noindent {\bf Corollary 3.5.} Let $ \BB_n=\BB_P \oplus \BB_Q \,, $            
$$
\multline
\text{where}\quad \BB_P= End_{FS_n}(W_P)= 
\bigoplus_{\lambda \in H_0(k,l;n) } \BB_{\lambda} ,\\ 
\text{and}\quad 
\BB_Q =End_{FS_n}(W_Q)= \bigoplus_{\lambda \in 
H_1(k,l;n)} \BB_{\lambda}\\
(\text{and where}\; \BB_{\lambda}=End_{FS_n}(W_{\lambda} )\,),
\endmultline
$$
Similarly, 
$$
\AA _n=\Cal A _P \oplus \Cal A _Q\,, 
$$
$$
\text{where}\quad \Cal A _P =End_{FA_n}(W_P)\;\text{and}\;
\Cal A _Q =End_{FA_n}(W_Q)                              
$$
Then $\AA _Q= \BB _Q\,.$
\bigskip
\noindent {\bf Proof.} Recall from (3.2.1) that 
$\varphi ^* ( FS_n)= P \oplus Q\,$ and $PW_Q=0\,.$ Therefore
$$
\BB _Q=End_{P\oplus Q}(W_Q)=End_{Q}(W_Q)\,.
$$
Similarly, $\varphi ^* (FA_n)=P'\oplus Q'$ and $P'W_Q=0\,,$ while by Lemma 3.4,
$Q'=Q\,.$ Thus also
$$
\AA _Q=End_{P'\oplus Q'}(W_Q)=End_{Q'}(W_Q)=End_{Q}(W_Q)\ \,.
$$
\smallskip
\noindent {\bf Note} that Corollary 3.5 proves part (b) of Theorem 3.2.

\bigskip
\bigskip
\noindent  {\bf \S 4 Asymptotics, as $n\to\infty\,.$}
\bigskip
\noindent We examine first the case $\dim V_0\ne \dim V_1\,.$ Let 
$$
P,\;\;Q,\;\; P',\;\; Q',\;\;\AA _n,\;\; \AA _P,\;\; \AA _Q,\;\; 
\BB _n,\;\;\BB _P,\;\; 
\BB _Q,\;\; W_P \quad\text{and}\quad W_Q
$$
as in Section 3.
\bigskip
\noindent {\bf Remark 4.1.} Let $\dim V_0\ne \dim V_1\,,$
say, $ \dim V_1 < \dim V_0\,.$ Then
$$
\lim_{n\longrightarrow \infty}
\frac{\dim P}{\dim\varphi ^* ( FS_n)}=0\,\qquad
\text{and}\qquad
\lim_{n\longrightarrow \infty}
\frac{\dim Q}{\dim\varphi ^* ( FS_n)}=1\,, \tag 4.1.1
$$
and similarly
$$
\lim_{n\longrightarrow \infty}
\frac{\dim P'}{\dim\varphi ^* ( FA_n)}=0\,\qquad
\text{and}\qquad
\lim_{n\longrightarrow \infty}
\frac{\dim Q'}{\dim\varphi ^* ( FA_n)}=1\,. \tag 4.1.2
$$
Together with $Q=Q',$ this clearly implies 
$$
\lim_{n\longrightarrow \infty}
\frac{\dim\varphi ^* ( FA_n)}{\dim\varphi ^*  ( FS_n)}= 1\,. \tag 4.1.3
$$
\medskip
\noindent {\bf The proof} follows since, as $n \to \infty\,,$ it can be
shown (by arguments similar to those in [BR \S 7]) that the 
growth of $\dim P$ is close to $(2\dim V_1)^n\,,$ while that of
$\dim Q$ is close to $(\dim V_0+\dim V_1)^n\,.$ This implies (4.1.1). The
proof of (4.1.2) follows by similar arguments. We skip the details.

\noindent It is also possible to show that here,

$$
\lim_{n \longrightarrow \infty}
\frac{\dim \BB_n}{
\dim \Cal A _n} =1         \tag 4.1.4
$$
and also

$$
\lim_{n \longrightarrow \infty}
\frac{\dim W_P}{
\dim V^{\otimes n}} =0\,\qquad
\text{and}\qquad
\lim_{n \longrightarrow \infty}
\frac{\dim W_Q}{
\dim V^{\otimes n}} =1\,. \tag 4.1.5
$$
Again, we skip the details.
\bigskip
\noindent {\bf 4.2.} In the rest of this section we study the case 
$\dim V_0=\dim V_1=k.$
\medskip
\noindent It follows from (2.4.1) and (2.6.3) that
$$
\multline
\varphi ^*  (FS_n)=M_1\oplus M_2 ,\qquad\text{where}\\
M_1\cong\left [\bigoplus_{\lambda '< \lambda \in H(k,k;n) }
(I_{\lambda} \oplus I_{\lambda '})\right ]\quad \text{and}\quad 
M_2\cong\left [\bigoplus_{\lambda '=\lambda \in H(k,k;n)}
I_{\lambda} \right ].                                      
\endmultline\tag 4.2.1
$$
\noindent {\bf Denote} 
$$
\multline
H_{sc}(k,k;n)=\{\lambda\in H(k,k;n)\; | \; \lambda=\lambda ' \}
\quad  (\text{`$sc$' for `self-conjugate')},\\ 
\text{and}\quad H_{nsc}(k,k;n)=H(k,k;n)\backslash H_s(k,k;n)
\quad  (\text{`$nsc$' for `non-self-conjugate'}).
\endmultline
$$
Thus, (4.2.1) can be rewritten with 
$$
M_1\cong\left [\bigoplus_{\lambda \in H_{nsc}(k,k;n) }
I_{\lambda}\right ]\quad \text{and}\quad 
M_2\cong\left [\bigoplus_{\lambda \in H_{sc}(k,k;n)}
I_{\lambda} \right ].   
$$
We compare $\dim M_1$ with $\dim M_2$, as well as 
the cardinalities of 
$H_{sc}(k,k;n)$ and of $H_{nsc}(k,k,;n),$ as $n\to\infty\,.$
Here we have
\bigskip
\noindent {\bf Proposition 4.3.} As $n \longrightarrow \infty$, 
$$
\frac{|H_{sc}(k,k;n)|}{|H_{nsc}(k,k;n)|}\approx
\frac{|H_{sc}(k,k;n)|}{|H(k,k;n)|}\approx
c_1\cdot \left (\frac{1}{n}\right ) ^k\quad \text{where}
\quad c_1= \frac{k\cdot (2k-1)!}{2^{k-1}}\,,\tag 4.3.1)
$$
and
$$
\frac{\dim M_2}{\dim M_1} \approx 
\frac{\dim M_2}{\dim \varphi ^* (FS_n)} \approx c_2 \cdot 
\left (\frac{1}{\sqrt n }\right ) ^k\,,\tag 4.3.2
$$
where $c_2$ is a constant that can be calculated explicitly by
Selberg integrals.
\smallskip
\noindent In particular, 
$$
\frac{\dim M_2}{\dim \varphi ^* (FS_n)} \approx 0\quad\text{and}\quad
\frac{\dim M_1}{\dim \varphi ^* (FS_n)} \approx 1. \tag 4.3.3
$$
\smallskip

\noindent {\bf Proof.} First, let 
$H'_{sc}(k,k;n)=\{(\lambda _1,\lambda _2,\cdots )
\in H_{sc}(k,k;n) \;|\; 
 \lambda _k \ge k\}$ (i.e $k\times k \subseteq \lambda ).$
It can be shown, for example by [BR 7.3] and 
by induction on $k,$ that $|H_{sc}(k,k;n)|\approx |H'_{sc}(k,k;n)|.$ 
It is easy to see, for example by the arguments below, that
$|H(k,k;n)|\approx |H_{nsc}(k,k;n)|.$ Therefore, to prove 
(4.3.1) it suffices to show that 
$$
\frac{|H'_{sc}(k,k;n)|}{|H(k,k;n)|}\approx
c_1\cdot \left (\frac{1}{n} \right )^k. \tag 4.3.1'
$$
It follows from [BR 7.3] that 
$$
|H(k,k;n)|\approx \frac{n^{2k-1}}{k!k!(2k -1)!}\,.
$$
Let $Par_k(m)$ denote the number of the partitions of $m$ to at most
$k$ parts. It is easy to see that 
$$
|H'_{sc}(k,k;n)|=Par_k((n-k^2)/2)\,.
$$
Let $m=\frac{n-k^2}{2},$ then $m\approx \frac{n}{2}.$ Again by  [BR 7.3]
(with $l=0$), $Par_k(m)\approx \frac{m^{k-1}}{k!(k-1)!}\,.$
Thus 
$$
H'_{sc}(k,k;n)|\approx \frac{n^{k-1}}{k!(k-1)!2^{k-1}}\,,
$$
and the proof of (4.3.1') follows.
%
%
\medskip
To prove (4.3.2), follow the rest of Section 7 in [BR]. 
Let $\lambda ' = \lambda \in H(k,k;n).$
Since $\lambda ' = \lambda ,$ in the notations of [BR 7.9 -- 7.21] ,
$a_i=b_i$ and $\alpha_i=\beta_i,\; i=1,\cdots, k.$
By 7.16.1 there, 
$$
d_{\lambda}\approx c(k,k)\cdot D_k^2(\alpha_1, \cdots, \alpha_k) \cdot
e^{-2k(\sum \alpha_i^2)}\cdot \left (\frac{1}{n} \right )
^{\theta (k,k)}\cdot (2k)^n, \tag 4.3.4
$$
where  $\theta (k,k)=\frac{1}{2}(k(k+1)-1).$

\noindent As in the proof of [BR, Thm 7.18], we follow [R,\S 2] 
to deduce that 
$$
\dim M_2 = \sum_{\lambda '=\lambda \in H(k,k;n)} d_{\lambda}^2\approx
\left [c(k,k)\cdot \left(\frac{1}{n}\right ) 
^{\theta (k,k)}\cdot (2k)^n\right ]^2 \cdot
(\sqrt n)^{k-1} \cdot I(k,k,2), \tag 4.3.5
$$
where 
$$
I(k,k,2)=\int_{P(k)} \left [ D_k^2(x_1,\cdots,x_k) \cdot
e^{-2k(\sum x_i^2)}
\right ]^2 d^{(k-1)}(x) \tag 4.3.6
$$
with
$$
P(k)=\{(x_1, \cdots , x_k) \, | \, 
x_1 \ge \cdots \ge x_k \; \text{and} \; \sum x_i=0 \}
$$
(correction: the term $\sum \alpha_i^2+ \sum \beta _j^2$ 
just  below (7.18.1) in [BR] --  should have been 
$\sum x_i^2+ \sum y_j^2$).
On the other hand, [BR, Thm 7.18] gives the asymptotics of 
$S_{k,k}^{(2)}(n)=\dim M:$

as $ n \longrightarrow \infty, $
$$
\dim M \approx \left [ c(k,k) \cdot \left ( \frac{1}{n} \right )
^\theta(k,k) \cdot (2k)^n \right ] ^2 \cdot
(\sqrt n)^{2k-1} \cdot I'(k,k,2) \tag 4.3.7
$$
for a certain multi-integral $I'(k,k,2)$ (given explicitly in [BR (7.18.1)]).

\noindent The proof of  (4.3.2) now 
clearly follows from  (4.3.5) and (4.3.7). \hfill \qed
\medskip

\bigskip
\noindent {\bf Remark 4.4.} Recall that 
$$
\multline
M_1=M_1(k,n)  \cong\left [\bigoplus_{\lambda '< \lambda \in H(k,k;n) }
(I_{\lambda} \oplus I_{\lambda '})\right ]\quad \text{and}\\ 
M_2=M_2(k,n) \cong\left [\bigoplus_{\lambda '=\lambda \in H(k,k;n)}
I_{\lambda} \right ].                                      
\endmultline
$$
Similarly, by removing the restriction that
$\lambda \in H(k,k;n),$ one can define 

$$
M_1(n) =\left [\bigoplus_{\lambda '< \lambda \vdash n }
(I_{\lambda} \oplus I_{\lambda '})\right ]\quad \text{and}\quad 
M_2(n)=\left [\bigoplus_{\lambda '=\lambda \vdash n }
I_{\lambda} \right ]. \tag 4.4.1                                      
$$
In analogy with Proposition 4.3, it should be interesting to calculate the 
asymptotics of the ratios
$$
\frac{\dim M_2(n)}{\dim M_1(n)} \qquad \text{and} 
\qquad \frac{\dim M_2(n)}{n!}
$$ 
as $n \longrightarrow \infty,$ {\it if} such asymptotics exists -- which 
is not at all obvious.

%
  
\bigskip
\bigskip
\centerline{\bf References}
\bigskip
 
\noindent [BR]  Berele A., Regev A.,~~Hook Young diagrams with
applications to combinatorics and to representations of Lie
superalgebras, {\it Adv. Math.} {\bf 64} (1987) 118--175.
\medskip
\noindent [B] Boerner H., Representations of groups, North-Holland, 
Amsterdam, 1963.
\medskip
\noindent [J] James, G. D.,  The representations of the
symmetric groups, Springer LNM, No. 682, Springer-Verlag, 1978.
\medskip

\noindent [JK] James, G. D. and Kerber A., The representation theory of the
symmetric group, Encyclopedia of Mathematics and its Applications, Vol. 16,
Addison-Wesley, 1981.
\medskip
\noindent [M] Macdonald I. G., Symmetric functions and Hall polynomials, 
2nd edition, Oxford University Press, 1995. 
\medskip

\noindent [R] ~~ Regev A., ~~Asymptotic values for degrees associated 
with strips
of Young diagrams,

\noindent {\it Adv. Math.} {\bf 41} (1981) 115--136.
\medskip
\noindent [Sa] Sagan B. E., The symmetric group, representations, 
combinatorial algorithms, and symmetric functions, Wadsworth \& 
Brooks/Cole Mathematics Series, 1991.
\medskip
\noindent [Sc] $\quad$ Schur I., Uber die rationalen Darstellungen der 
allgemeinen linearen Gruppe (1927), {\it in} ``I.Schur, Gesammetle
Abhandlungen III,'' pp. 68-85, Springer-Verlag, Berlin, 1973.
\medskip
\noindent [Se] $\quad$ Sergeev A., The tensor algebra of the identity 
representation as a module over the Lie superalgebras $gl(n,m)$
and $Q(n),$  (Russian) Tom 123 (165) (1984). English translation:
Math, USSR Sbornik Vol. 51 (1985), No. 2, 419-427.
\medskip
\noindent [W] Weyl H., The Classical groups, Princeton Math. Series, No. 1, 
Princeton Univ. Press, Princeton, N.J.

\bye

\bye